\documentclass[a4paper,11pt,reqno]{amsart}

\usepackage{amsmath,amssymb,amsthm,mathrsfs}
\usepackage{txfonts}

\usepackage{tikz-cd}
\usetikzlibrary{calc,intersections,through,backgrounds}
\usetikzlibrary{decorations.markings}
\usetikzlibrary{arrows.meta}

\usepackage{hyperref}

\makeatletter
\@namedef{subjclassname@2020}{%
  \textup{2020} Mathematics Subject Classification}
\makeatother

\numberwithin{equation}{section}
\usepackage[scr=boondoxo]{mathalfa}
\newcommand{\cA}{\mathscr{A}}
\newcommand{\cB}{\mathscr{B}}
\newcommand{\cF}{\mathscr{F}}
\newcommand{\cX}{\mathscr{X}}

\newtheorem{Thm}{Theorem}[section]
\newtheorem{Cor}[Thm]{Corollary}
\newtheorem{Lem}[Thm]{Lemma}
\newtheorem{Pro}[Thm]{Proposition}
\theoremstyle{definition}
\newtheorem{Def}[Thm]{Definition}

\newcommand{\N}{\mathbb{N}}
\newcommand{\R}{\mathbb{R}}
\newcommand{\Z}{\mathbb{Z}}

\newcommand{\al}{\alpha}
\newcommand{\eps}{\varepsilon}
\newcommand{\del}{\delta}
\newcommand{\Del}{\Delta}
\newcommand{\gam}{\gamma}
\newcommand{\lam}{\lambda}
\renewcommand{\rho}{\varrho}

\newcommand{\om}{\omega}

\newcommand{\asrk}{\operatorname{asrk}}
\newcommand{\comb}{\text{\rm comb}}

\newcommand{\diam}{\operatorname{diam}}
\newcommand{\es}{\emptyset}
\newcommand{\id}{\operatorname{id}}
\renewcommand{\d}{\partial}
\newcommand{\sm}{\setminus}
\newcommand{\sub}{\subset}

\newcommand{\hh}{hierarchically hyperbolic}
\renewcommand{\SS}{\operatorname{SS}}
\newcommand{\GL}{\operatorname{GL}}
\newcommand{\SL}{\operatorname{SL}}
\newcommand{\Or}{\operatorname{O}}
\newcommand{\li}{l_{\infty}}
\newcommand{\compl}{\text{\rm c}}
\begin{document}

\title{A combinatorial higher-rank hyperbolicity condition}
\author{Martina J{\o}rgensen}
\address{Department of Mathematics\\ETH Z\"urich\\R\"amistrasse 101
  \\8092 Z\"urich\\Switzerland}
\email{martina.joergensen@math.ethz.ch}
\author{Urs Lang}
\address{Department of Mathematics\\ETH Z\"urich\\R\"amistrasse 101
  \\8092 Z\"urich\\Switzerland}
\email{lang@math.ethz.ch}
\date{30 September 2023}
\thanks{Research supported by Swiss National Science Foundation Grant 197090.}

\begin{abstract}
We investigate a coarse version of a $2(n+1)$-point inequality 
characterizing metric spaces of combinatorial dimension at most~$n$
due to Dress. This condition, experimentally called $(n,\del)$-hyperbolicity,
reduces to Gromov's quadruple definition of $\del$-hyperbolicity in
case $n = 1$. The $l_\infty$-product of $n$ $\del$-hyperbolic spaces
is $(n,\del)$-hyperbolic. Every $(n,\del)$-hyperbolic metric space,
without any further assumptions, possesses a slim $(n+1)$-simplex property
analogous to the slimness of quasi-geodesic triangles in Gromov hyperbolic
spaces. In connection with recent work in geometric group theory, we show
that every Helly group and every \hh\ group of (asymptotic) rank $n$ acts
geometrically on some $(n,\del)$-hyperbolic space.
\end{abstract}

\maketitle


\section{Introduction}

Generalizations and variations of Gromov hyperbolicity~\cite{G87}
belong to the most present themes in geometric group theory today
(see, for example, the introduction in~\cite{HKS1} for a comprehensive
list of these developments). Here we continue the investigation
of higher-rank hyperbolicity phenomena from \cite{G93}~(Sect.~6.B$_2$),
\cite{W11}, \cite{KL}, and~\cite{GL}. These results show
in particular that most of the characteristic properties of
Gromov hyperbolic spaces, regarding the shape of triangles,
quasi-geodesics, and isoperimetric inequalities, among others,
have adequate rank~$n$ analogues ($n \ge 2$) in a context of generalized
global non-positive curvature. We refer to the paragraph preceding
Theorem~\ref{Thm:ss} below for a sample result.
The focus in the present paper is on a more foundational, partly
combinatorial aspect. We explore a coarse \mbox{$2(n+1)$-point} inequality
for general metric spaces that reduces to Gromov's quadruple definition of
$\del$-hyperbolicity in the case $n = 1$ and, if $\del = 0$, to an
inequality characterizing metric spaces of combinatorial dimension
at most~$n$ due to Dress~\cite{Dr}. The latter concept measures the
combinatorial complexity of the induced metric on finite subsets in terms
of the dimension of the (polyhedral) injective hull of these sets
(see below and Sect.~\ref{sect:ih}). Throughout the paper, the unified
condition is referred to briefly as {\em $(n,\del)$-hyperbolicity}.
This notion turns out to possess a variety of remarkable
properties, tying up higher-rank hyperbolicity with (coarsely) injective
metric spaces, injective hulls, and some recent developments in geometric
group theory.

We now proceed to the actual definition. The common origin of the two
aforementioned special cases $\del = 0$ and $n = 1$ is the well-known
observation that a metric space $X$ admits an isometric embedding into a
metric ($\R$-)tree if and only if every quadruple $(x,x',y,y')$ of points
in $X$ satisfies the inequality
\begin{equation} \label{eq:tree}
  d(x,x') + d(y,y')
  \le \max \bigl\{ d(x,y) + d(x',y'), d(x,y') + d(x',y) \bigr\}.
\end{equation}
The smallest complete such tree is provided by the
{\em injective hull} $E(X)$~\cite{I} of $X$; see pp.~322, 326, and~329
in~\cite{Dr}, where the injective hull is referred to as the {\em tight span}
$T_X$ of $X$. A metric space $Y$ is {\em injective}, as an object in the
metric category with $1$-Lipschitz maps as morphisms, if every such map
$\phi \colon A \to Y$ defined on a subset of a metric space $B$ extends
to a $1$-Lipschitz map $\bar\phi \colon B \to Y$. The injective hull $E(X)$
is characterized by the universal property that every isometric embedding
of $X$ into some injective metric space $Y$ factors through $E(X)$.
On the one hand, Theorem~9 on p.~327 in~\cite{Dr} generalizes the above
observation to metric spaces with an injective hull of dimension at most $n$
or, more precisely, with the property that the topological dimension of $E(V)$
is less than or equal to $n$ for every finite set $V \sub X$. The respective
$2(n+1)$-point condition is precisely what we call $(n,0)$-hyperbolicity
in Definition~\ref{def:nhyp} below (and we will give another proof of Dress'
theorem in Corollary~\ref{Cor:dress}). On the other hand, Definition~1.1.C
in~\cite{G87} is equivalent to the relaxed inequality~\eqref{eq:tree}:
$X$ is {\em $\del$-hyperbolic}, for $\del \ge 0$, if and only if 
\begin{equation} \label{eq:hyp}
  d(x,x') + d(y,y')
  \le \max \bigl\{ d(x,y) + d(x',y'), d(x,y') + d(x',y) \bigr\} + 2\del
\end{equation}
for every quadruple $(x,x',y,y')$ of points in $X$. In the case $n = 1$, the
following condition is indeed a (somewhat inefficient) reformulation
of this inequality (see Proposition~\ref{Pro:hyp} for the details).

\begin{Def} \label{def:nhyp}
Let $n \ge 0$ be an integer, and let
$I = I_n$ denote the $2(n+1)$-element set $\{\pm 1,\pm 2,\ldots,\pm(n+1)\}$
with the canonical involution $-\id$.
A metric space $X$ is called {\em $(n,\del)$-hyperbolic}, for some
$\del \ge 0$, if for every family $(x_i)_{i \in I}$ of points in $X$
there exists a permutation $\al \ne -\id$ of $I$ such that
\begin{equation} \label{eq:nhyp}
  \sum_{i\in I}d(x_i,x_{-i}) \le \sum_{i\in I}d(x_i,x_{\alpha(i)}) + 2\del.
\end{equation}
We say that $X$ is {\em $(n,\ast)$-hyperbolic} if $X$ is $(n,\del)$-hyperbolic
for some $\del$.
\end{Def}

To emphasize the analogy with~\eqref{eq:hyp}, we could choose
$\alpha$ in~\eqref{eq:nhyp} so as to maximize the sum on the right.
Note, however, that for $n = 1$ there are $|I| = 4$ summands on either side.
If $n = 0$, then $\al = \id$ is the only permutation of
$I = \{1,-1\}$ distinct from $-\id$, thus a metric space $X$ is
$(0,\del)$-hyperbolic if and only if the diameter $\diam(X)$ is less than
or equal to $\del$. An $(n,\del)$-hyperbolic metric space is
$(n',\del')$-hyperbolic for all $n ' \ge n$ and $\del' \ge \del$
(see Lemma~\ref{Lem:mon}).

We briefly summarize some further basic properties. 
The $\li$-product of $(n_i,\del)$-hyperbolic
spaces $X_i$, $i = 1,2$, is $(n_1+n_2,\del)$-hyperbolic. In particular,
the $\li$-product of $n$ $\del$-hyperbolic spaces is
$(n,\del)$-hyperbolic (Proposition~\ref{Pro:prod}).
In general, $(n,\ast)$-hyperbolicity is preserved under rough isometries,
that is, $(1,c)$-quasi-isometries for $c \ge 0$ (Lemma~\ref{Lem:ri}).
Anticipating Theorem~\ref{Thm:cinj}, we mention that quasi-isometry invariance
is granted for the class of coarsely injective metric spaces.
(It should be noted that for general metric spaces, quasi-isometry invariance
fails also in the case $n = 1$; see, for example, Remark 4.1.3(2) in~\cite{BS}.)
The asymptotic rank $\asrk(X)$ of an $(n,\del)$-hyperbolic
space $X$ is at most $n$, and every asymptotic cone of a sequence of pointed
$(n,\del)$-hyperbolic spaces is $(n,0)$-hyperbolic
(Proposition~\ref{Pro:asrk}).

Next, we relate the notion of $(n,\del)$-hyperbolicity to injective hulls.
Remarkably, Gromov's $\del$-inequality~\eqref{eq:hyp} passes on from $X$
to $E(X)$.
This provides a most efficient way of embedding a general $\del$-hyperbolic
metric space into a complete, contractible, geodesic $\del$-hyperbolic space
with some more features reminiscent of global non-positive curvature
(see Sect.~4.4 in~\cite{DMT} and Propositions~1.2, 1.3, and~3.8 in~\cite{L}).
Likewise, the injective hull of an $(n,\del)$-hyperbolic
metric space $X$ is $(n,\del)$-hyperbolic (Proposition~\ref{Pro:ih}).
We then prove the following characterization. The key step is the implication
$(4) \Rightarrow (1)$, which is shown in a quantitative form in
Proposition~\ref{Pro:copy}.

\begin{Thm} \label{Thm:nhyp-ih}
For every metric space $X$ and every integer $n \ge 0$,
the following assertions are equivalent:
\begin{enumerate}
\item
$X$ is $(n,\ast)$-hyperbolic;
\item
the injective hull $E(X)$ is $(n,\ast)$-hyperbolic;
\item
$\asrk(E(X)) \le n$;
\item
there is a constant $r_0$ such that $E(X)$ contains no isometric copy
of the \mbox{$\li$-ball}
$B(0,r) = [-r,r]^{n+1} \sub \li^{\,n+1} := (\R^{n+1}, \|\cdot\|_\infty)$
for $r > r_0$.
\end{enumerate}
\end{Thm}

We turn to a more geometric higher-rank hyperbolicity condition, analogous
to the slimness of quasi-geodesic triangles in Gromov hyperbolic spaces.
We say that a metric space $X$ has the {\em slim simplex property $(\SS_n)$}
if for all $\lam \ge 1$ and $c \ge 0$ there exists a constant $D \ge 0$
such that if $\Del$ is a Euclidean
$(n+1)$-simplex and $\phi \colon \d\Del \to X$ is a map whose restriction
to every facet is a $(\lam,c)$-quasi-isometric embedding, then the image
of every facet is within distance $D$ of the union of the images of the
remaining ones. This property was first established in Theorem~1.1 in~\cite{KL}
for spaces of asymptotic rank at most~$n$ in a setup reminiscent of 
non-positive curvature, including in particular all proper metric spaces
with a conical geodesic bicombing (as defined in~\cite{DL15}).
A stronger uniform statement has been given in Theorem~7.2 in~\cite{GL},
and by virtue of the properties of the injective hull we can deduce a
completely general result in the present context.

\begin{Thm} \label{Thm:ss}
Every $(n,\del)$-hyperbolic metric space $X$ satisfies the slim simplex
property $(\SS_n)$ with a constant $D$ depending only on $n,\del,\lam,c$.
\end{Thm}

In fact, the argument yields a constant of the form
$D = (1+c) \cdot D'(n,\del,\lam)$ (see Theorem~\ref{Thm:sss}).
Proposition~7.4 in~\cite{GL} shows in turn
that every metric space $X$ satisfying~$(\SS_n)$ with
$D = (1+c) \cdot D'(\lam)$ has asymptotic rank at most $n$.
However, it is not true in general that a metric space satisfying 
$(\SS_n)$ is $(n,\ast)$-hyperbolic. For example, $l_\infty^{\,n}$ is
$(n,0)$-hyperbolic and hence has property $(\SS_n)$,
whereas the Euclidean $\R^n$ (being quasi-isometric to $l_\infty^{\,n}$)
satisfies $(\SS_n)$ but fails to be $(n,\ast)$-hyperbolic for $n \ge 2$
(see Proposition~\ref{Pro:norms}). On the positive side, it follows easily
from the implication $(4) \Rightarrow (1)$ in Theorem~\ref{Thm:nhyp-ih} that
every injective metric space $X$ with property $(\SS_n)$ is
$(n,\ast)$-hyperbolic.
Since $(n,\ast)$-hyperbolicity is preserved under rough isometries,
it is thus natural to seek a generalization to the following class of
metric spaces, recently considered in~\cite{CCGHO},
\cite{H1}, and~\cite{HHP}. We call a metric space $X$ {\em coarsely injective}
if there is a constant $c \ge 0$ such that every 1-Lipschitz map
$\phi \colon A \to X$ defined on a subset of a metric space $B$ has
an extension $\bar \phi \colon B \to X$ satisfying
$d(\bar\phi(b),\bar\phi(b')) \le d(b,b') + c$ for all $b,b' \in B$.
This is equivalent to $E(X)$ being within finite distance of the image
of the canonical embedding $e \colon X \to E(X)$ and also to $X$ being
roughly isometric to an injective metric space; see
Proposition~\ref{Pro:c-inj}. It was shown in~\cite{L99}
(see also~\cite{CE}, \cite{L}) that every geodesic Gromov hyperbolic space
is coarsely injective. Thus the following result generalizes various 
known characterizations of hyperbolicity to higher rank.

\begin{Thm} \label{Thm:cinj}
Let $X$ be a coarsely injective metric space. For every $n \ge 0$,
the following properties are equivalent:
\begin{enumerate}
\item
$X$ is $(n,\ast)$-hyperbolic;
\item
$X$ satisfies the slim simplex property $(\SS_n)$;
\item
$\asrk(X) \le n$; 
\item
every asymptotic cone of\/ $X$ is $(n,0)$-hyperbolic;  
\item
for all $c > 0$ there exists $r_0$ such that there is no
$(1,c)$-quasi-isometric embedding of\/ $B(0,r) \sub \li^{\,n+1}$
into $X$ for $r > r_0$.
\end{enumerate}
If\/ $X$ is in addition proper and cocompact, then 
the asymptotic rank of\/ $X$ is finite, so all properties hold for
$n = \asrk(X)$.
\end{Thm}

Note that since conditions~(2) and~(3) are preserved under quasi-isometries,
the theorem also shows that $(n,\ast)$-hyperbolicity is a quasi-isometry
invariant for coarsely injective spaces. 

We now discuss a few applications of the above results in connection
with some recent developments in geometric group theory.

A first corollary pertains to the class of Helly groups introduced
in~\cite{CCGHO} and further explored in~\cite{H2}, \cite{HuO},
and \cite{OV}. A connected locally finite graph is called
a {\em Helly graph} if its vertex set $V$, endowed with the natural integer
valued metric, has the property that every family of pairwise
intersecting balls has non-empty intersection. Then $V$ is coarsely injective,
and the injective hull $E(V)$ is proper and has the structure of a locally
finite polyhedral complex with only finitely many isometry types of $n$-cells,
isometric to injective polytopes in $l_\infty^{\,n}$, for every $n \ge 1$.
Furthermore, if the graph has uniformly bounded degrees, then $E(V)$
is finite-dimensional. See (the proofs of) Theorem~1.1 in~\cite{L} and
Theorem~6.3 in~\cite{CCGHO}. 
A group $G$ is called a {\em Helly group} if $G$ acts geometrically
(that is, properly discontinuously and cocompactly by isometries)
on the vertex set $V$ of a Helly graph and, hence, geometrically 
on $E(V)$. The following corollary of Theorems~\ref{Thm:nhyp-ih}
and~\ref{Thm:ss} applies more generally to groups acting geometrically
on an injective metric space. Recent examples in~\cite{HV} show that
not every such group is Helly (see Corollary~D therein).
Recall also that a group with a geometric action on a proper
geodesic metric space $X$ is finitely generated and, with respect to any word
metric, quasi-isometric to $X$ (see, for example, Theorem~8.37 in~\cite{DK}).

\begin{Cor} \label{Cor:helly}
Suppose that $G$ is a group acting geometrically on a proper injective
metric space $X$, and endow $G$ with any word metric.
Then $n := \asrk(X) = \asrk(G)$ is finite and agrees with
\begin{enumerate}
\item[(1)] the minimal integer $n_1$ such that $X$ is $(n_1,\ast)$-hyperbolic;
\item[(2)] the minimal integer $n_2$ such that $G$
  satisfies~$(\SS_{n_2})$;
\item[(3)] the maximal integer $n_3$ such that there is a quasi-isometric
  embedding of $\R^{n_3}$ into $G$;
\item[(4)] the maximal integer $n_4$ such that 
  $X$ contains an isometric copy of\/~$\li^{\,n_4}$.
\end{enumerate}
Furthermore, $G$ has no free abelian subgroup of rank $n+1$.
\end{Cor}

Next we relate $(n,\ast)$-hyperbolicity to \hh\ spaces or groups, as defined
in~\cite{BHS17} and~\cite{BHS19}. Every \hh\ space $(X,d)$ has a finite
rank~$\nu$ (see Definition~1.10 in~\cite{BHS21}). In general,
\[
  \nu \le \bar\nu \le \asrk(X),
\]
where $\bar\nu$ denotes the supremum
of all integers $k$ such that there exist constants $\lam,c$ and
$(\lam,c)$-quasi-isometric embeddings of $B(0,r) \sub \R^k$ into $X$
for all $r > 0$ (a quasi-isometry invariant).
However, according to the discussion in Sect.~1.1.3 in~\cite{BHS21},
the equality $\nu = \bar\nu$ holds for all natural examples
of \hh\ spaces and in particular for all \hh\ groups. For the latter,
$\bar\nu$ equals the maximal dimension of a quasi-flat in the group.
In the recent paper~\cite{HHP}, Theorem~A, it is shown that every \hh\ space
$(X,d)$ admits a coarsely injective metric $\rho$ quasi-isometric to $d$.
In the case of a \hh\ group $X = G$, the metric $\rho$ can be chosen
so that $G$ acts geometrically on the (proper) coarsely injective space
$(X,\rho)$. Combining this result with Theorem~\ref{Thm:cinj}, we get
the following corollary.

\begin{Cor} \label{Cor:hh}
Let $(X,d)$ be a \hh\ group $X = G$ of rank $\nu$ or,
more generally, a \hh\ space with $\nu = \bar\nu$, and let $\rho$
be a coarsely injective metric on $X$ quasi-isometric to $d$ (see above).
Then $\nu = \asrk(X,\rho) = \asrk(X,d)$, and this is the minimal integer
such that $(X,\rho)$ is $(\nu,\ast)$-hyperbolic. In particular, $X$ satisfies
the slim simplex property~$(\SS_\nu)$ with respect to either $\rho$ or $d$.
\end{Cor}

Lastly, we consider Riemannian symmetric spaces of non-compact type.
Spaces of rank $1$ are Gromov hyperbolic and thus $(1,\ast)$-hyperbolic in the
above terminology. Spaces of rank $n \ge 2$ have asymptotic rank $n$ and
satisfy the slim simplex property $(\SS_n)$ by~\cite{KL}; however, with
respect to the Riemannian metric, they are not $(n,\ast)$-hyperbolic.
The only $n$-dimensional $(n,\ast)$-hyperbolic normed space is
$l_\infty^{\,n}$, up to isometry (see Proposition~\ref{Pro:norms}),
so the question is whether a given non-compact symmetric space $X = G/K$
of rank $n \ge 2$ admits an $(n,\ast)$-hyperbolic, $G$-invariant
distance function such that the maximal flats are isometric to $l_\infty^{\,n}$
with respect to the induced metric. A result in~\cite{PT} (see also~\cite{P})
shows that the $G$-invariant distance functions on~$X$ corresponding to norms
on the maximal flats are in bijection with the $G$-invariant Finsler metrics
(of class $C^0$) on $X$ and also with the norms on $T_{\!p}F$ invariant under
the Weyl group, for a base point $p$ of $X$ and any maximal flat $F$
through $p$. Thus the rank $n$ symmetric spaces in question are those
whose Weyl group preserves an $n$-cube, and it remains to see that the
resulting metric is indeed $(n,\ast)$-hyperbolic. 
The recent paper~\cite{H1} shows that for $X = \GL(n,\R)/\Or(n)$,
as well as for every classical irreducible symmetric space
of non-compact type associated with the automorphism group $G$ of a
non-degenerate bilinear or sesquilinear form, there is a coarsely injective,
$G$-invariant metric $d$ on $X$ making the maximal flats isometric
to $l_\infty^{\,n}$; furthermore, the injective hull of $(X,d)$ is proper.
The following immediate consequence of Theorem~\ref{Thm:cinj} thus applies
to all classical groups $G$ {\em not of type $\SL$},
as defined in~\cite{H1}.

\begin{Cor} \label{Cor:symm}
If a symmetric space $X = G/K$ of non-compact type and rank $n \ge 2$
is equipped with a $G$-invariant coarsely injective Finsler
metric $d$, then $(X,d)$ is $(n,\ast)$-hyperbolic.
\end{Cor}

The rest of the paper is organized as follows.
Sect.~\ref{sect:basics} records the basic properties of $(n,\ast)$-hyperbolic
spaces. In Sect.~\ref{sect:ih} we first review the construction of the
injective hull and the definition of the combinatorial dimension,
then we prove some auxiliary results. Sect.~\ref{sect:ih-nhyp}
discusses injective hulls of $(n,\ast)$-hyperbolic spaces and establishes
Theorem~\ref{Thm:nhyp-ih}. In the concluding Sect.~\ref{sect:ss-ci}
we turn to the slim simplex property and prove the remaining result stated
above.


\section{Basic properties} \label{sect:basics}

We begin with some elementary observations regarding Definition~\ref{def:nhyp}.
First we notice that in case $n \ge 1$ the permutation $\alpha \ne -\id$
can always be taken to be fixed point free. Recall that we put
$I_n = \{\pm 1,\ldots,\pm(n+1)\}$.

\begin{Lem} \label{Lem:fixed-pts}
Let $X$ be a metric space, and let $x_i \in X$ for $i \in I = I_n$, where
$n \ge 1$. Then for every permutation $\alpha \ne -\id$ of $I$ there
is a permutation $\bar\al \ne -\id$ of $I$ without fixed points such that
$S(\al) := \sum_{i \in I} d(x_i,x_{\al(i)}) \le S(\bar\al)$.
\end{Lem}

\begin{proof}
This holds trivially for $\al = \id$, as $S(\id) = 0$.
On the other hand, if $\al \ne \id$
and $\al(i) = i$ for some $i \in I$, then $\al$ has a cycle involving a pair
$j \ne k$ with $\al(j) = k$, and there is an $\al' \ne -\id$ 
with $\al'(j) = i$ and $\al'(i) = k$ such that $S(\al') \ge S(\al)$ by the
triangle inequality. Eliminating all fixed points in this way, one gets a
permutation $\bar\al$ as desired.
\end{proof}

We now give the details of the characterization for $n = 1$.

\begin{Pro} \label{Pro:hyp}
A metric space $X$ is $(1,\del)$-hyperbolic if and only if\/ $X$ is Gromov
$\del$-hyperbolic. 
\end{Pro}

\begin{proof}
If $X$ is $\del$-hyperbolic, then by adding the term
$L := d(x,x') + d(y,y')$ to each of the three sums in~\eqref{eq:hyp}
and by substituting $(x_1,x_{-1},x_2,x_{-2}) := (x,x',y,y')$, one
gets~\eqref{eq:nhyp} for some cyclic permutation $\al \ne -\id$ of
$I = \{\pm 1,\pm 2\}$.

Conversely, suppose that~\eqref{eq:nhyp} holds for some $\al \ne -\id$.
By Lemma~\ref{Lem:fixed-pts} we can assume that $\al$ has no fixed points.
We use the identification $(x,x',y,y') := (x_1,x_{-1},x_2,x_{-2})$
and consider three cases. In the first case, $\al$ is cyclic, and the
respective sum $S(\al)$ equals $L$ plus either $d(x,y) + d(x',y')$
or $d(x,y') + d(x',y)$. Then~\eqref{eq:hyp} follows upon subtracting $L$
on both sides. In the second case, $\al$ is still cyclic, but
$S(\al) = d(x,y) + d(y,x') + d(x',y') + d(y',x)$.
Then there is an involution $\hat\al \ne -\id$ such that
$S(\al) \le S(\hat\al)$, the latter sum being equal to
$2(d(x,y) + d(x',y'))$ or $2(d(x,y') + d(x',y))$.
This reduces the second case to the remaining case, where $\al$
is an involution distinct from $-\id$. Then, dividing~\eqref{eq:nhyp} by~$2$,
one obtains~\eqref{eq:hyp} with $\del$ in place of $2\del$.
See Figure~\ref{fig:perm} for illustration.
\begin{figure}
\smallskip
\begin{tikzpicture}[scale=1]
\draw[line width=1pt,black!35!white] (0,0) rectangle (1,1);  
\fill (0,0) circle (1.35pt) node[anchor=east]{$x_{-1}$};
\fill (1,0) circle (1.35pt) node[anchor=west]{$x_1$};
\fill (0,1) circle (1.35pt) node[anchor=east]{$x_{-2}$};
\fill (1,1) circle (1.35pt) node[anchor=west]{$x_2$};
\draw[line width=1pt,black!35!white] (2,0)--(3,0)--(2,1)--(3,1)--cycle;
\fill (2,0) circle (1.35pt);
\fill (3,0) circle (1.35pt);
\fill (2,1) circle (1.35pt);
\fill (3,1) circle (1.35pt);
\draw[line width=1pt,black!35!white] (4,0)--(4,1)--(5,0)--(5,1)--cycle;
\fill (4,0) circle (1.35pt);
\fill (4,1) circle (1.35pt);
\fill (5,0) circle (1.35pt);
\fill (5,1) circle (1.35pt);
\draw[line width=1pt,black!35!white] (6,0)..controls (5.8,0.5)..(6,1)..controls (6.2,0.5)..cycle;
\draw[line width=1pt,black!35!white] (7,0)..controls (6.8,0.5)..(7,1)..controls (7.2,0.5)..cycle;
\fill (6,0) circle (1.35pt);
\fill (6,1) circle (1.35pt);
\fill (7,0) circle (1.35pt);
\fill (7,1) circle (1.35pt);  
\draw[line width=1pt,black!35!white] (8,0)..controls (8.33,0.67)..(9,1)..controls (8.67,0.33)..cycle;
\draw[line width=1pt,black!35!white] (9,0)..controls (8.33,0.33)..(8,1)..controls (8.67,0.67)..cycle;
\fill (8,0) circle (1.35pt);
\fill (8,1) circle (1.35pt);
\fill (9,0) circle (1.35pt);
\fill (9,1) circle (1.35pt);
\end{tikzpicture} 
\caption{Illustration of the sums $S(\al)$ for the fixed point free
  permutations $\al \ne -\id$ of $\{\pm 1,\pm 2\}$.}
\label{fig:perm}
\end{figure}
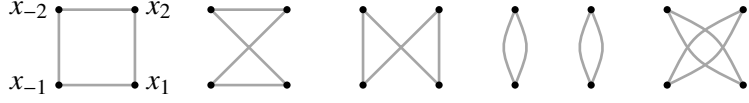
\end{proof}

For completeness we record the obvious monotonicity properties.

\begin{Lem} \label{Lem:mon}
Let $X$ be a metric space.
\begin{enumerate}
\item
If\/ $X$ is $(n,\del)$-hyperbolic, then $X$ is $(n',\del)$-hyperbolic for all
$n' > n$;
\item
$X$ is $(n,\del)$-hyperbolic if and only if\/ $X$ is $(n,\del')$-hyperbolic
for all $\del' > \del$.
\end{enumerate}
\end{Lem}  

\begin{proof}
For~(1), given $n' > n$ and a family of points $x_i$ with $i \in I_{n'}$,
there is a permutation $\al \ne -\id$ of $I_n$ such that~\eqref{eq:nhyp} holds
for the corresponding subfamily, and we can simply add the
terms $d(x_i,x_{-i})$ for $i \in I_{n'} \sm I_n$ on both sides and extend
$\al$ accordingly.

To prove the `if' direction in~(2),
let $x_i \in X$ for $i \in I = I_n$. For every $\del' > \del$ there is a
permutation $\al \ne -\id$ of~$I$ such that~\eqref{eq:nhyp} holds with
$\del'$ in place of $\del$, and for a suitable sequence $\del'_k \to \del$
the corresponding permutations all agree, so that~\eqref{eq:nhyp} holds
in the limit.
\end{proof}  

We turn to products.
The {\em $\li$-product} of an $l$-tuple of metric spaces $(X_k,d_k)$,
$k = 1,\ldots,l$, is the set $X = \prod_{k=1}^l X_k$ with the metric $d$
defined by
\[
d(x,y) := \max \bigl\{ d_k(x_k,y_k): k = 1,\ldots,l \bigr\}
\]
for all pairs of points $x = (x_1,\ldots,x_l)$ and $y = (y_1,\ldots,y_l)$
in $X$. The following proposition is a direct adaptation of the
result for $\del = 0$ given in~\cite{Dr}, (5.15). 

\begin{Pro} \label{Pro:prod}
If\/ $(X,d)$ is the $\li$-product of an $l$-tuple of metric spaces
as above, and if\/ $(X_k,d_k)$ is $(n_k,\del)$-hyperbolic,
then $(X,d)$ is $(n,\del)$-hyperbolic for $n := \sum_{k=1}^l n_k$.
In particular, the $\li$-product of $n$ $\del$-hyperbolic metric spaces
is $(n,\del)$-hyperbolic.
\end{Pro}

\begin{proof}
Let $x_i = (x_{i,1},\ldots,x_{i,l}) \in X$ for $i \in I = I_n$. 
For $k = 1,\ldots,l$, define
\[
I(k) := \bigl\{ i \in I: d(x_i,x_{-i}) = d_k(x_{i,k},x_{-i,k}) \bigr\}.
\]
Note that $I = \bigcup_{k=1}^l I(k)$ and $I(k) = -I(k)$, in particular
$|I(k)|$ is even. Since $|I| = 2(n+1) > 2n$, there is an index $k$ with
$|I(k)| > 2n_k$ and hence $|I(k)| \ge 2(n_k + 1)$.
As $(X_k,d_k)$ is $(n_k,\del)$-hyperbolic, there exists a permutation
$\al \ne -\id$ of $I(k)$ such that 
\[
\sum_{i\in I(k)} d_k(x_{i,k},x_{-i,k})
\le \sum_{i\in I(k)} d_k(x_{i,k},x_{\al(i),k}) + 2\del
\]
(compare the first part of Lemma~\ref{Lem:mon} in case $|I(k)| > 2(n_k + 1)$).
Using the definition of $I(k)$ and the inequality $d_k(x_{i,k},x_{\al(i),k})
\le d(x_i,x_{\al(i)})$, and extending $\al$ by $-\id$ on $I \sm I(k)$,
we get a permutation $\al \neq -\id$ of $I$ such that~\eqref{eq:nhyp} holds.
\end{proof}

Let $X = (X,d)$ and $Y = (Y,d)$ be two metric spaces.
For constants $\lam,c \ge 0$, a map $\phi \colon X \to Y$
will be called {\em $(\lam,c)$-Lipschitz} if
\[
  d(\phi(x),\phi(x')) \le \lam\,d(x,x') + c
\]
for all $x,x' \in X$. If, in addition, $\lam \ge 1$ and
\[
  d(\phi(x),\phi(x')) \ge \lam^{-1}d(x,x') - c
\]
for all $x,x' \in X$, then $\phi$ is a
{\em $(\lam,c)$-quasi-isometric embedding}; in the case $\lam = 1$
we call $\phi$ a {\em roughly isometric embedding}.
A {\em quasi-isometry} or {\em rough isometry}
$\phi \colon X \to Y$ is a quasi-isometric or roughly isometric embedding,
respectively, such that $Y$ is within finite distance from the
image $\phi(X)$. 

\begin{Lem} \label{Lem:ri}
If $Y$ is an $(n,\del)$-hyperbolic space and $f \colon X \to Y$ is a
$(1,\eps)$-quasi-isometric embedding,
then $X$ is $(n,\del + 2(n+1)\eps)$-hyperbolic.
In particular, $(n,\ast)$-hyperbolicity is preserved under rough isometries.
\end{Lem}

\begin{proof}
Let $x_i \in X$ for $i \in I = I_n$. Then
\begin{align*}
\sum_{i\in I}d(x_i,x_{-i})
&\leq \sum_{i\in I} d(f(x_i),f(x_{-i})) + 2(n+1)\eps\\
&\leq\sum_{i\in I} d(f(x_i),f(x_{\al(i)})) + 2\del + 2(n+1)\eps\\
&\leq\sum_{i\in I} d(x_i,x_{\al(i)}) + 2\del + 4(n+1)\eps
\end{align*}
for some permutation $\al \ne -\id$ of $I$.
\end{proof}

Evidently every $(n,\ast)$-hyperbolic normed space (and, more generally,
every metric space admitting dilations) is in fact
$(n,0)$-hyperbolic. The following classification follows from some
well-known results about injective hulls of normed spaces together
with Dress' characterization of the combinatorial dimension, but
can also be proved more directly.

\begin{Pro} \label{Pro:norms}
A normed space is $(n,0)$-hyperbolic if and only if it
is finite-dimensional with a polyhedral norm, in which case the minimal
such $n$ equals the number of pairs of opposite facets of the unit ball.
In particular, every $(n,0)$-hyperbolic normed space has dimension at most
$n$, and equality occurs if and only if it is
isometric to $\li^{\,n}$.
\end{Pro}

\begin{proof}
Let $(X,\|\cdot\|)$ be a finite-dimensional normed space with a polyhedral
norm, such that the unit ball $B$ has $n$ pairs
$\pm F_1,\ldots,\pm F_n$ of opposite facets.
For each of these pairs, let $f_i \in X^*$ be the linear functional that
is $\pm 1$ on $\pm F_i$. Then $f = (f_1,\ldots,f_n) \colon X \to l_\infty^{\,n}$
is a linear isometric embedding, and since $l_\infty^{\,n}$ is $(n,0)$-hyperbolic
by Proposition~\ref{Pro:prod}, so is $X$. To see that $n$ is minimal with this
property, choose a relatively interior point $x_i$ in each $F_i$, and put
$x_{-i} := -x_i$. This gives a set $\{x_i: i \in I_{n-1}\} \sub \d B$ of
cardinality $2n$ such that no two distinct elements are connected by a
line segment in $\d B$.
Then $\|x_i - x_j\| = 2\,\bigl\| \frac12(x_i + x_{-j}) \bigr\| < 2$
whenever $j \ne -i$, and so
\[
  \sum_{i \in I_{n-1}}\|x_i - x_{\al(i)}\| < \sum_{i \in I_{n-1}}\|x_i - x_{-i}\|
\]
for every permutation $\al \ne -\id$ of $I_{n-1}$. Thus $X$ is not
$(n-1,0)$-hyperbolic. Clearly $n$ is greater than or equal to the dimension
of $X$, with equality if and only if $B$ is a parallelotope and $X$ is
isometric to $l_\infty^{\,n}$ via the above $f$.

Suppose now that the unit ball $B$ is not polyhedral, whereas, without loss
of generality, $X$ is still finite-dimensional.
Choose a convex set $C_1 \sub \d B$ that is maximal
with respect to inclusion (a singleton if $B$ is strictly convex),
and a point $x_1$ in the interior $C_1^0$ relative
to the affine hull of $C_1$. Then no point in $\d B \sm C_1$ is connected
to $x_1$ by a line segment in $\d B$. Recursively, for $n \ge 2$, if $C_{n-1}$
and $x_{n-1} \in C_{n-1}^0$ are chosen, pick a maximal convex set $C_n \sub \d B$
such that $C_n^0$ is disjoint from
$D_{n-1} = \bigcup_{k=1}^{n-1} (C_k \cup -C_k)$, and a point $x_n \in C_n^0$.
Note that $D_{n-1} \neq \es$ for all $n$, because $B$ is not polyhedral.
Thus, for arbitrarily large $n$, we find a set
$\{\pm x_1,\ldots,\pm x_n\} \sub \d B$ such that no two distinct elements
are connected by a line segment in $\d B$, and it follows as above
that $X$ is not $(n-1,0)$-hyperbolic.
\end{proof}

Next we relate $(n,\ast)$-hyperbolicity to the notion of asymptotic rank,
which originates from~\cite{G93} and was further discussed
in~\cite{K}, \cite{W11}, \cite{D}.

Given a sequence $(X_k)_{k\in\N}$ of metric spaces $X_k = (X_k,d_k)$,
we call a compact metric space $Z = (Z,d_Z)$ an {\em asymptotic subset} of
$(X_k)$ if there exist a sequence $0 < r_k \to \infty$ and subsets
$Z_k \sub X_k$ such that the rescaled sets $(Z_k,r_k^{-1}d_k)$
converge in the Gromov--Hausdorff topology to $Z$; equivalently,
there exist sequences of positive numbers $r_k \to \infty$, $\eps_k \to 0$,
and $(1,\eps_k)$-quasi-isometric embeddings
$\phi_k \colon Z \to (X_k,r_k^{-1}d_k)$.
Every asymptotic subset admits an isometric embedding
into some asymptotic cone $X_\om$ of $(X_k)$ with the same scale
factors (where $\om$ is any non-principal ultrafilter on $\N$) and, conversely,
every compact subset of an asymptotic cone $X_\om$ of $(X_k)$
is an asymptotic subset of some subsequence of $(X_k)$
(see Sect.~10.6 in~\cite{DK} for a discussion of asymptotic cones).
We define the {\em asymptotic rank} of the sequence $(X_k)$ as the supremum
of all integers $m \ge 0$ for which there exists an asymptotic subset of
$(X_k)$ bi-Lipschitz homeomorphic to a compact subset of $\R^m$ with
positive Lebesgue measure. It can be shown by a metric differentiation
argument that if such an asymptotic subset exists, then there is a norm on
$\R^m$ whose unit ball is an asymptotic subset of some subsequence of
$(X_k)$ (see Corollary~2.2 and Proposition~3.1 in~\cite{W11}).
For a single metric space $X = (X,d)$, the {\em asymptotic rank\/} $\asrk(X)$
of $X$ is defined as the asymptotic rank of the constant sequence
$(X_k,d_k) = (X,d)$.

\begin{Pro} \label{Pro:asrk}
Let $(X_k)_{k \in \N}$ be a sequence of\/ $(n,\del)$-hyperbolic metric spaces
$X_k = (X_k,d_k)$. Then every asymptotic subset and every asymptotic cone
of\/ $(X_k)$ is $(n,0)$-hyperbolic, and the asymptotic rank
of\/ $(X_k)$ is at most~$n$. In particular, $\asrk(X) \le n$ for any
$(n,\ast)$-hyperbolic space $X$.
\end{Pro}

\begin{proof}
Let $Z$ be an asymptotic subset of $(X_k)$.
There are sequences $r_k \to \infty$ and $\eps_k \to 0$ such that, for every
$k$, there exists a $(1,\eps_k)$-quasi-isometric embedding
of $Z$ into the $(n,r_k^{-1}\del)$-hyperbolic space 
$(X_k,r_k^{-1}d_k)$.
Thus $Z$ is $(n,r_k^{-1}\del + 2(n+1)\eps_k)$-hyperbolic
for all $k$ (Lemma~\ref{Lem:ri}) and hence $(n,0)$-hyperbolic
(Lemma~\ref{Lem:mon}).

If $X_\om$ is an asymptotic cone of $(X_k)$,
then every finite set $Z \sub X_\om$ is an asymptotic subset of some
subsequence of $(X_k)$, hence $Z$ is $(n,0)$-hyperbolic, and so is $X_\om$. 

For the assertions about the asymptotic rank, suppose that $Z$ is an
asymptotic subset of $(X_k)$ bi-Lipschitz homeomorphic to a compact subset
of $\R^m$ with positive Lebesgue measure. Then, as mentioned above, there is a
norm $\|\cdot\|$ on $\R^m$ whose unit ball $B$ is an asymptotic subset of
some subsequence of $(X_k)$. By the first part of the proof, $B$ is
$(n,0)$-hyperbolic, and so $m \le n$ by Proposition~\ref{Pro:norms}.
This shows that the asymptotic rank of $(X_k)$ is at most $n$.
\end{proof}


\section{Injective hulls and combinatorial dimension} \label{sect:ih}

In this section we first review the definitions of the injective hull
and the combinatorial dimension, then we state some auxiliary results.

Recall that a metric space $(Y,d)$ is {\em injective} if partially defined
$1$-Lipschitz maps into $Y$ can always be extended while preserving the
Lipschitz constant.
The most basic examples of injective metric spaces are $\R$, $\li(S)$
for any index set~$S$, complete $\R$-trees, and $\li$-products thereof.
Injective metric spaces are complete, geodesic, contractible, and share
some more properties with spaces of non-positive curvature (see~\cite{L}).
Isbell~\cite{I} showed that every metric space has an
{\em injective hull\/} $(e,E(X))$, that is, 
$E(X)$ is an injective metric space, $e \colon X \to E(X)$ is an isometric
embedding, and for every isometric embedding of $X$ into some injective metric
space $Y$ there is an isometric embedding $E(X) \to Y$ such that the
diagram
\[
\begin{tikzcd}[arrows={-stealth}]
  X \rar["e"] \drar & E(X) \dar \\
  & Y
\end{tikzcd}
\]
commutes.
If $(i,Y)$ is another injective hull of $X$, then there exists a unique
isometry $j \colon E(X) \to Y$ such that $j \circ e = i$.
Isbell's construction was rediscovered and further investigated by
Dress~\cite{Dr} who called $E(X)$ the {\em tight span} $T_X$ of $X$.
We briefly review the explicit construction of $E(X)$.

Let $\Del(X)$ denote the set of all functions $f \colon X \to \R$ satisfying
\[
  f(x) + f(y) \ge d(x,y)
\]
for all $x,y \in X$ (in~\cite{Dr}, $\Del(X)$ is denoted $P_X$). The subset
$E(X) \sub \Del(X)$ of {\em extremal functions} (in the terminology of~\cite{I})
consists of all minimal elements of the partially ordered set $(\Del(X),\le)$.
If $f \in \Del(X)$, then
\begin{equation} \label{eq:f*}
  f^*(x) := \sup_{z \in X}\,(d(x,z) - f(z)) \le f(x) 
\end{equation}
for all $x \in X$, and $q(f) := \frac12(f + f^*) \le f$ belongs to $\Del(X)$.
Hence, if $f \in E(X)$, then $f^* = f$, and, conversely, every function
$f \colon X \to \R$ with $f^* = f$ is extremal. Moreover, by iterating
the transformation $q \colon \Del(X) \to \Del(X)$ and by passing to pointwise
limits, one obtains a canonical map
\[
  p \colon \Del(X) \to E(X)
\]
such that $p(f) \le f$ for all $f \in \Del(X)$ and
$p(f \circ \gam) = p(f) \circ \gam$ for all isometries $\gam$ of $X$
(see p.~332 in~\cite{Dr} or Proposition~3.1 in~\cite{L}).
For every $y \in X$, the distance function $d_y := d(\cdot,y)$ is extremal.
By plugging the inequality $d(x,z) \le d(x,y) + d(y,z)$
into the definition in~\eqref{eq:f*} one sees that $f^*$ is $1$-Lipschitz,
thus every $f \in E(X)$ is $1$-Lipschitz and satisfies  
$d_y - f(y) \le f \le d_y + f(y)$ and   
\begin{equation} \label{eq:f-d}
  \|f - d_y\|_\infty = \sup_{x \in X}|f(x) - d_y(x)| = f(y) 
\end{equation}
for all $y \in X$.
It follows that $\|f - g\|_\infty < \infty$ for every pair of functions
$f,g \in E(X)$. This provides $E(X)$ with a metric, and the map
\[
  e \colon X \to E(X), \quad x \mapsto d_x,
\]
is a canonical isometric embedding.
The retraction $p \colon \Del(X) \to E(X)$ is $1$-Lipschitz
with respect to the (possibly infinite) $\li$-distance on $\Del(X)$.
We refer to~\cite{I} and~\cite{L} for two different proofs
that $(e,E(X))$ is indeed an injective hull of $X$.

Let $f,g \in E(X)$. It follows from~\eqref{eq:f-d} and the
triangle inequality that
\begin{equation} \label{eq:x-f-g-y}
  f(x) + \|f-g\|_\infty + g(y) \ge \|d_x - d_y\|_\infty = d(x,y)
\end{equation}
for all $x,y \in X$. The next lemma shows that the points $x,y$ can be
chosen such that equality holds up to an arbitrarily small error
(compare Theorem~3(iii) in~\cite{Dr}). In particular, if $X$ is compact,
then there exists a pair $x,y$ such that $f,g$ lie on a
geodesic from $d_x$ to $d_y$.

\begin{Lem} \label{Lem:xfgy}
For $f,g \in E(X)$ and $\eps > 0$, there exist $x,y \in X$ such that
\[
  f(x) + \|f-g\|_\infty + g(y) < d(x,y) + \eps.
\]
\end{Lem}  

\begin{proof}
Pick $x,y$ such that one of the following two conditions holds:
$\|f-g\|_{\infty} < f(y) - g(y) + \frac{\eps}{2}$
and $f(y) < d(x,y) - f(x) + \frac{\eps}{2}$, or
$\|f-g\|_{\infty} < g(x) - f(x) + \frac{\eps}{2}$ and
$g(x) < d(x,y) - g(y) + \frac{\eps}{2}$. In either case, this
gives the desired inequality.  
\end{proof}

We remark further that if $X$ is compact, then so is $E(X)$,
as a direct consequence of the Arzel\`a--Ascoli theorem.

Suppose now, for the moment, that $X$ is finite.
Then $E(X)$ is a subcomplex of the boundary of the unbounded polyhedral
set $\Del(X) \sub \R^X$. The faces of $\Del(X)$ that belong to
$E(X)$ are exactly those whose affine hull is determined by a system
of equations of the form
\[
  f(x_i) + f(x_j) = d(x_i,x_j)
\]
such that every point of $X$ occurs at least once in the system.
These are precisely the bounded faces of $\Del(X)$.
Given $f \in E(X)$, there is a unique minimal face~$P$ containing
$f$ in its relative interior. The dimension of $P$ can be read off from the
{\em equality graph} of $f$ on the set $X$ with edges $\{x_i,x_j\}$
corresponding to the above equations:
$f$ is uniquely determined on all connected components with a cycle of
odd length or a loop $\{x,x\}$ (occurring only if $f = d_x$), whereas $f$
has one degree of freedom on each of the remaining components. Thus 
$\dim(P) = n$ is the number of bipartite connected components.
If $x_1,\ldots,x_n \in X$ are such that there
is one $x_i$ in each of them, then the map sending $g \in P$ to
$(g(x_1),\ldots,g(x_n))$ is an isometry from $P$ onto a polytope in
$\li^{\,n}$ (see Lemma~4.1 and Theorem~4.3 in~\cite{L}).
In particular, for finite $X$, $E(X)$ has the structure of a finite polyhedral
complex of dimension at most $\frac12|X|$ with $\li$-metrics on
the cells. See p.~93 in~\cite{Dr+} for the possible shapes of $E(X)$ for
generic metric spaces up to cardinality $5$.

The {\em combinatorial dimension} $\dim_\comb(X)$ of a metric space $X$
equals the supremum of $\dim(E(V))$ over all finite subsets $V \sub X$.
This notion was introduced by Dress in~\cite{Dr}. Theorem~$9'$ (on p.~380)
therein provides a variety of characterizations, whereas Theorem~9
in the introduction highlights the equivalence of $X$ being
$(n,0)$-hyperbolic, in our terminology, and the inequality
$\dim_\comb(X) \le n$. This equivalence will also follow from the results
in the next section in combination with the following characterization
(see Corollary~\ref{Cor:dress}).

\begin{Pro} \label{Pro:dim-comb}
For every metric space $X$ and $n \ge 1$, the following are equivalent:
\begin{enumerate}
\item   
  $\dim_\comb(X) \ge n$;
\item
  $E(X)$ contains an isometric copy of a non-empty open subset
  of\/ $\li^{\,n}$;
\item
  $E(X)$ contains an isometric copy of\/
  $\{\pm se_i: i = 1,\ldots,n\} \sub \li^{\,n} = (\R^n,\|\cdot\|_\infty)$
  for some $s > 0$ (where $e_1,\ldots,e_n$ denote the canonical basis vectors
  of\/ $\R^n$ as usual).
\end{enumerate}
\end{Pro}

\begin{proof}
If (1) holds, then there is a finite set $V \sub X$ with
$\dim(E(V)) \ge n$, and $E(V)$ embeds isometrically into   
$E(X)$, thus $E(X)$ contains a copy of an $n$-dimensional
cell of $E(V)$, and (2)~follows. Evidently (2) implies~(3).

We show that~(3) implies (1). Let $I := \{\pm 1,\ldots,\pm n\}$,
and let $\{f_i: i \in I\}$ be a subset of $E(X)$ isometric to
$\{\pm se_i: i = 1,\ldots,n\} \sub \li^{\,n}$ for some $s > 0$,
so that $\|f_i - f_{-i}\|_\infty = 2s$ for all $i \in I$ and
$\|f_i - f_j\|_\infty = s$ whenever $j \in I \sm \{i,-i\}$.
Fix an $\eps \in (0,\frac{s}{2})$.
Lemma~\ref{Lem:xfgy} together with~\eqref{eq:x-f-g-y} shows that for every
$i \in \{1,\dots,n\}$ there exist $x_i,x_{-i} \in X$ such that
\[
  d(x_i,x_{-i}) \le f_i(x_i) + 2s + f_{-i}(x_{-i}) < d(x_i,x_{-i}) + \eps.
\]
Note that $d(x_i,x_{-i}) > 2s - \eps > 0$.
Furthermore, if $j \in I \sm \{i,-i\}$, then
\[
  d(x_i,x_j) \le f_i(x_i) + s + f_j(x_j);
\]
assuming without loss of generality that $f_i(x_i) \ge f_j(x_j)$, we
infer that
\begin{align*}
  d(x_i,x_j)
  &\ge d(x_i,x_{-i}) - d(x_j,x_{-i}) \\
  &> \bigl( f_i(x_i) + 2s + f_{-i}(x_{-i}) - \eps \bigr)
    - \bigl( f_j(x_j) + s + f_{-i}(x_{-i}) \bigr) \\
  &\ge s - \eps > 0.
\end{align*}
Thus the set $V := \{x_i: i \in I\}$ has cardinality $2n$. By putting
$h(x_i) := f_i(x_i) + s$ for all $i \in I$, we get a function $h \in \Del(V)$.
Let $g := p(h) \in E(V)$, and recall that $g \le h$.
For every $i \in I$, we have
$h(x_i) + h(x_{-i}) - \eps < d(x_i,x_{-i}) \le g(x_i) + g(x_{-i})$, hence
$g(x_i) > h(x_i) - \eps$, and if $j \in I \sm \{i,-i\}$, then
\[
  g(x_i) + g(x_j) > h(x_i) + h(x_j) - 2\eps \ge d(x_i,x_j) + s - 2\eps
  > d(x_i,x_j)
\]
by the choice of $\eps$. Since $g \in E(V)$, it follows that
$g(x_i) + g(x_{-i}) = d(x_i,x_{-i})$ for $i = 1,\ldots,n$.
Thus the equality graph of $g$ on $V$ has just $n$ pairwise disjoint edges,
hence $n$ bipartite components, and so the minimal cell of $E(V)$ containing
$g$ has dimension~$n$, as discussed earlier.
Since $\dim(E(V)) \le \frac12|V| = n$, we have $\dim(E(V)) = n$.
\end{proof}

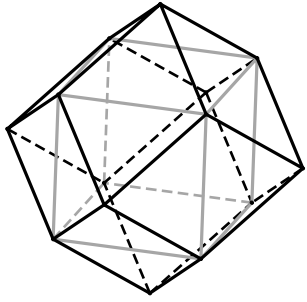
\begin{figure}
\begin{tikzpicture}[scale=0.45,xslant=-0.1]
	
\draw[xslant=0.1, white] (-8,-5) rectangle (5,5);	

\coordinate (A) at (-3.23,-1.09, 0);
\coordinate (B) at (1.03,-1.67, 0);
\coordinate (C) at (1.61, 2.59, 0);
\coordinate (D) at (-2.65, 3.17, 0);
\coordinate (A') at (-3.23, -1.09, 4.3);
\coordinate (B') at (1.03, -1.67, 4.3);
\coordinate (C') at (1.61, 2.59, 4.3);
\coordinate (D') at (-2.65, 3.17, 4.3);	
	
\draw[line width=1.2pt,black!35!white] (A')--(B') (B')--(B) (D')--(C') (C')--(C) (B')--(C') 
(D')--(D) (C)--(D) (A')--(D') (B)--(C);

\begin{scope}[on background layer]
\draw[line width=1pt, black!35!white, dash pattern=on 4pt off 2pt, opacity=0.7] 
(A')--(A) (B)--(A) (A)--(D);
\end{scope}							
	
\coordinate (E) at (-5.07, 1.33, 2.15);
\coordinate (F) at (-0.81, 0.75, 6.45);
\coordinate (G) at (-1.39, -3.51, 2.15);
\coordinate (H) at (-0.23, 5.01, 2.15);
\coordinate (I) at (3.45, 0.17, 2.15);
\coordinate (J) at (-0.81, 0.75, -2.15);
	
\draw[line width=1.2pt,black] (F)--(A') (F)--(B') (F)--(C') (F)--(D')
(H)--(C') (H)--(D') (I)--(B') (I)--(C') (G)--(A') (G)--(B')
(I)--(C) (H)--(C) (H)--(D) (E)--(A') (E)--(D') (E)--(D);
		
\begin{scope}[on background layer]							  
\draw[line width=1pt, black, dash pattern=on 4pt off 2pt, opacity=0.7]
(J)--(A) (J)--(B) (J)--(C) (J)--(D) (G)--(A) (G)--(B) (I)--(B) (E)--(A);
\end{scope}		  								  
	
\filldraw [black] 
(A) circle (1.35pt)
(B) circle (1.35pt)
(C) circle (1.35pt)
(D) circle (1.35pt)
(A') circle (1.35pt)
(B') circle (1.35pt)
(C') circle (1.35pt)
(D') circle (1.35pt)
(E) circle (1.35pt)
(F) circle (1.35pt)
(G) circle (1.35pt)
(H) circle (1.35pt)
(I) circle (1.35pt)
(J) circle (1.35pt);
	
\end{tikzpicture}
\caption{The injective hull of $\{\pm se_i\} \sub \li^{\,3}$ is a rhombic
  dodecahedron whose extra vertices span a cube of edge length~$s$.}
\label{fig:rd}
\end{figure}

We conclude this section with a  quantitative version of the
above implication $(3) \Rightarrow (2)$.
If an injective metric space $Y$ contains an isometric copy of a set
$\{\pm se_i\} \sub \li^{\,n}$ as in~(3), then $Y$ also contains
an isometric copy of the injective hull $E(\{\pm se_i\})$. The latter is
isometric to a convex polytope in $\li^{\,n}$ which in turn contains the
ball $B\bigl(0,\frac s2\bigr) = \bigl[-\frac s2 ,\frac s2\bigr]^n
\sub \li^{\,n}$. In case $n = 3$, this polytope is the
rhombic dodecahedron shown in Figure~\ref{fig:rd} (compare~\cite{GM}).
For our purposes, the following result will suffice.

\begin{Lem} \label{Lem:ball}
Let $Y$ be an injective metric space. If\/ $Y$ contains an isometric copy
of\/ $\{\pm se_i\} \sub \li^{\,n}$ for some $s > 0$ and $n \ge 1$,
then there is also an isometric embedding of the ball
$B\bigl(0,\frac s2\bigr) = \bigl[-\frac s2 ,\frac s2\bigr]^n
\sub \li^{\,n}$ into $Y$.
\end{Lem}

\begin{proof}
Since $Y$ is injective, every isometric embedding
$\phi \colon \{\pm se_i\} \to Y$ extends to a $1$-Lipschitz map
$\phi \colon \{\pm se_i\} \cup \bigl[-\frac s2,\frac s2 \bigr]^n \to Y$.
Let $x,y \in \bigl[-\frac s2,\frac s2 \bigr]^n$.
After possibly interchanging $x$ and $y$, we have
$\|x - y\|_\infty = x_i - y_i$ for some $i$.
Then 
\begin{align*}
d(\phi(se_i),\phi(-se_i))
&\le d(\phi(se_i),\phi(x)) + d(\phi(x),\phi(y)) + d(\phi(y),\phi(-se_i)) \\
&\le \|se_i - x\|_\infty + \|x - y\|_\infty + \|y + se_i\|_\infty \\
&= (s - x_i) + (x_i - y_i) + (y_i + s),
\end{align*}
and since $d(\phi(se_i),\phi(-se_i)) = 2s$, the equality
$d(\phi(x),\phi(y)) = \|x - y\|_\infty$ holds.
\end{proof}


\section{The injective hull of an $(n,\del)$-hyperbolic space}
\label{sect:ih-nhyp}

We now relate $(n,\del)$-hyperbolicity to injective hulls.
The following proposition generalizes the known result for $n = 1$
(compare Sect.~4.4 in~\cite{DMT}, Chap.~5 in~\cite{Dr+},
or Proposition~1.3 in~\cite{L}). 

\begin{Pro} \label{Pro:ih}
A metric space $X$ is $(n,\del)$-hyperbolic if and only if its
injective hull $E(X)$ is $(n,\del)$-hyperbolic.
\end{Pro}

\begin{proof}
If $E(X)$ is $(n,\del)$-hyperbolic, then so is $e(X) \sub E(X)$
and hence $X$.

Conversely, suppose that $X$ is $(n,\del)$-hyperbolic, and let
$f_i \in E(X)$ for $i \in I = I_n$. Fix $\eps > 0$ for the moment.
Lemma~\ref{Lem:xfgy} shows that for every $i \in \{1,\dots,n+1\}$
there exist $x_i,x_{-i} \in X$ such that
\[
  \|f_i-f_{-i}\|_{\infty} < d(x_i,x_{-i}) - f_i(x_i) - f_{-i}(x_{-i}) + \eps.
\]
Setting $S := \frac12\sum_{i\in I}(f_i(x_i) + f_{-i}(x_{-i}))
= \sum_{i\in I}f_i(x_i)$, we get that
\[
  \sum_{i\in I}\|f_i-f_{-i}\|_{\infty} < \sum_{i\in I} d(x_i,x_{-i}) - 2S
  + 2(n+1)\eps.
\]
By the assumption on $X$ there exists a permutation $\al \ne -\id$ of $I$ with
\[
  \sum_{i\in I} d(x_i,x_{-i}) \le \sum_{i\in I}d(x_i,x_{\al(i)}) + 2\del.
\]
Using the inequalities $d(x_i,x_{\al(i)}) \le f_i(x_i) + f_i(x_{\al(i)})$ we
deduce that
\begin{align*}
\sum_{i\in I}\|f_i-f_{-i}\|_{\infty}
&< \sum_{i \in I}\bigl(f_i(x_i) + f_i(x_{\al(i)})\bigr) - 2S + 2\del + 2(n+1)\eps\\
&= \sum_{i \in I}\bigl(f_i(x_{\al(i)}) - f_{\al(i)}(x_{\al(i)})\bigr) + 2\del + 2(n+1)\eps\\
&\le \sum_{i \in I}\|f_i-f_{\al(i)}\|_{\infty} + 2\del + 2(n+1)\eps.
\end{align*}
Thus $E(X)$ is $(n,\del + (n+1)\eps)$-hyperbolic for all $\eps > 0$
and hence $(n,\del)$-hyperbolic (Lemma~\ref{Lem:mon}).
\end{proof}	

Our next goal is to show that a metric space $X$ is $(n,\ast)$-hyperbolic
if and only if, intuitively, its injective hull has no large
$(n+1)$-dimensional subsets. To measure the size, we will use the sets
$\{\pm se_i: i = 1,\ldots,n+1\} \sub \li^{\,n+1}$ for $s > 0$
(compare Lemma~\ref{Lem:ball}). In preparation for the actual result,
Proposition~\ref{Pro:copy} below, we state the following criterion.

\begin{Lem} \label{Lem:copy}
Let $n \ge 1$, and let $V = \{x_i: i \in I = I_n\}$ be a metric space of
cardinality $2(n+1)$. Let $\cA$ denote the set of the $n+1$ pairs
$\{x_i,x_{-i}\}$, and let $\cA^\compl$ denote the set of all pairs
$\{x,y\} \not\in \cA$ of two distinct points in $V$. Suppose that there is
a function $f \colon V \to \R$ such that $f(x_i) + f(x_{-i}) = d(x_i,x_{-i})$
for all pairs $\{x_i,x_{-i}\} \in \cA$ and
\[
  s := \min\bigl\{ f(x) + f(y) - d(x,y): \{x,y\} \in \cA^\compl \bigr\} > 0.
\]
Then $f \in E(V)$, and there is an isometric embedding of\/
$\{0\} \cup \{\pm se_i\} \sub \li^{\,n+1}$ into $E(V)$ mapping $0$ to $f$.
\end{Lem}

\begin{proof}
To show that $f \in E(V)$ it only remains to verify that $f \ge 0$.
For every $y \in V$ there is an edge $\{x_i,x_{-i}\} \in \cA$ not
containing $y$, so that $\{x_i,y\},\{x_{-i},y\} \in \cA^\compl$, hence
\begin{align*}
  f(x_i) + f(x_{-i}) &= d(x_i,x_{-i}) \\
                    &\le d(x_i,y) + d(x_{-i},y) \\
                    &\le f(x_i) + f(y) - s + f(x_{-i}) + f(y) - s
\end{align*}
and therefore $f(y) \ge s > 0$. Now, for every $i \in I$,
we define $f_i \colon V \to \R$ such that
$f_i(x_{\pm i}) = f(x_{\pm i}) \pm s$ and $f_i(x_j) = f(x_j)$ for all
$j \in I \sm \{i,-i\}$. Then $f_i(x_j) + f_i(x_{-j}) = d(x_j,x_{-j})$
for all $\{x_j,x_{-j}\} \in \cA$, and $f_i(x) + f_i(y) \ge d(x,y)$ whenever
$x,y \in V$; thus $f_i \in E(V)$.
Note that $\|f - f_i\|_\infty = s$ and $\|f_i - f_{-i}\|_\infty = 2s$
for all $i \in I$, and $\|f_i - f_j\|_\infty = s$ whenever $j \not\in \{i,-i\}$.
This yields an isometric embedding as desired.
\end{proof}

We now have the following key result.

\begin{Pro} \label{Pro:copy}		
Let $n \ge 1$. If\/ $X$ is $(n,\del)$-hyperbolic, then the injective hull
$E(X)$ contains no isometric copy of\/ $\{\pm se_i\} \sub \li^{\,n+1}$
for $s > \del$. Conversely, if\/ $E(X)$ contains no isometric copy of\/
$\{\pm se_i\} \sub \li^{\,n+1}$ for $s > \del$,
then $X$ is $(n,n\del)$-hyperbolic.
\end{Pro}

For $n = 1$, this reduces to the well-known fact that $X$ is $\del$-hyperbolic
if and only if $E(X)$ contains no isometric copy of
$\{\pm se_i\} \sub \li^{\,2}$, or, equivalently, of
$[0,s]^2 \sub \ell_1^{\,2}$, for $s > \del$ (compare p.~335f in~\cite{Dr},
the introduction in~\cite{BC}, and the discussion at the end of Sect.~3
in~\cite{DL16}).

\begin{proof}
Suppose that $X$ is $(n,\del)$-hyperbolic, and for some $s > 0$ there is a
subset $\{f_i: i \in I = I_n\}$ of $E(X)$ isometric to
$\{\pm se_i\} \sub \li^{\,n+1}$, so that $\|f_i - f_{-i}\|_\infty = 2s$
for all $i \in I$ and $\|f_i - f_j\|_\infty = s$ whenever $j \not\in \{i,-i\}$.
By Proposition~\ref{Pro:ih}, $E(X)$ is itself $(n,\del)$-hyperbolic, thus 
\[
  \sum_{i\in I}\|f_i - f_{-i}\|_\infty
  \le \sum_{i\in I}\|f_i - f_{\al(i)}\|_\infty + 2\del
\]
for some permutation $\al \ne -\id$ of $I$. Note that all summands in the
first sum are equal to the maximal distance $2s$, whereas at least
two terms in the second sum are $\le s$. It follows that $s \le \del$.

We prove the second part. If $|X| < 2(n+1)$, then $X$ is $(n,0)$-hyperbolic.
Suppose now that $V = \{x_i: i \in I = I_n\} \sub X$ is a set of cardinality
$2(n+1)$. Define $\cA$ and $\cA^\compl$ as in Lemma~\ref{Lem:copy},
and consider the set $\cF$ of all functions $f \colon V \to \R$ such
that $f(x_i) + f(x_{-i}) = d(x_i,x_{-i})$ for all $\{x_i,x_{-i}\} \in \cA$.
For $f \in \cF$, put
\begin{align*}
  s_f &:= \min\bigl\{ f(x) + f(y) - d(x,y): \{x,y\} \in \cA^\compl \bigr\},\\
  \cB_f &:= \bigl\{ \{x,y\} \in \cA^\compl: f(x) + f(y) - d(x,y) = s_f \bigr\}.
\end{align*}
Note that $s_f \le \diam(V)$, because there is a pair
$\{x_i,x_j\} \in \cA^\compl$ such that $f(x_i) \le \frac12d(x_i,x_{-i})$ and
$f(x_j) \le \frac12d(x_j,x_{-j})$. We now fix $f \in \cF$ for the rest of the
proof such that 
\[
  s_f = \bar s := \sup\{s_g: g \in \cF\}
\]
and $|\cB_f| \le |\cB_g|$ for all $g \in \cF$ with $s_g = \bar s$.
The elements of $\cA$ and $\cB := \cB_f$ will be called $\cA$-edges and
$\cB$-edges. We claim that for every $\cA$-edge, either both vertices belong
also to a $\cB$-edge, or neither of the two vertices has this property.
Suppose to the contrary that for some index $i \in I$, the point $x_i$ is
in $\bigcup\cB$, whereas $x_{-i}$ is not. Then, for some sufficiently small
$\eps > 0$, the function $g \in \cF$ defined by
$g(x_{\pm i}) := f(x_{\pm i}) \pm \eps$ and $g(y) := f(y)$ otherwise would
satisfy $s_g = \bar s$ and $\cB_g \subsetneqq \cB = \cB_f$, in contradiction to
the choice of~$f$.
Since $\cB \neq \es$, it follows from this claim that there is a non-empty
connected subgraph of $(V,\cA \cup \cB)$ such that each of its
vertices belongs to a unique $\cA$-edge and at least one $\cB$-edge.
Among all such subgraphs we select one with the least number of edges
and call it~$G$. There are two possible types, as described next.

The first possibility is that $G$ is simply a cycle graph with an even number
of edges alternating between $\cA$ and $\cB$. In this case we choose an
orientation of $G$ and define the permutation $\al \colon I \to I$ such that
the map $x_i \mapsto x_{\al(i)}$ sends each vertex of $G$ to the following one
and every other point $x_i$ in $V$ to $x_{-i}$.
In the remaining case, when $G$ is not an alternating cycle, by minimality
$G$ contains no such cycle as a proper subgraph either. Then, starting with
an oriented $\cA$-edge of $G$, we follow an alternating path in $G$,
stopping at the first vertex $v$ that was visited already earlier.
Since the $\cA$-edges are pairwise disjoint, the last edge belongs to $\cB$.
As there is no alternating cycle, by deleting the initial subpath up to the
first occurrence of $v$ we get an alternating loop based at $v$ that starts
and ends with a $\cB$-edge.
Proceeding with the oriented $\cA$-edge issuing from $v$,
we choose another alternating path, ending at the first vertex
$w$ occurring already earlier in the whole construction. Again, the last
edge is in $\cB$, and since there is no alternating cycle it follows that $w$
cannot be part of the loop based at $v$. We conclude that in the second
case, $G$ consists of two disjoint alternating loops based at $v$ and $w$,
respectively, each starting and ending with a $\cB$-edge, and an alternating
path from $v$ to $w$ (possibly of length one), starting and
ending with an $\cA$-edge. Then we define the permutation $\al \colon I \to I$
such that the map $x_i \mapsto x_{\al(i)}$ cyclically permutes each of the two
loops, moving every vertex to the next one, and interchanges the two
vertices of every $\cB$-edge in the path from $v$ to $w$.
Furthermore, as in the first case, $x_i \mapsto x_{\al(i)} = x_{-i}$
on the remaining part of~$V$. See Figure~\ref{fig:g}.
\begin{figure}
\begin{tikzpicture}[scale=0.73,>=Stealth]
		
\newcommand \x{1.5};
		
\draw[white] (-\x-0.5,-\x-0.5) rectangle (\x+0.5, \x+0.5);
		
\tikzset{->/.style={decoration={markings,mark=at position 1 with {\arrow[scale=0.83]{>}}},postaction={decorate}}}
		
\coordinate (A) at (\x,0);
\coordinate (B) at ({\x /2},{sqrt(3)*\x /2});
\coordinate (C) at ({-\x /2},{sqrt(3)*\x /2});
\coordinate (D) at (-\x,0);
\coordinate (E) at ({-\x /2},{-sqrt(3)*\x /2});
\coordinate (F) at ({\x /2},{-sqrt(3)*\x /2});
		
\foreach \a/\b in {A/B, C/D, E/F}
\draw[->,line width=1pt,black!35!white] (\a) -- (\b);	
		
\foreach \a/\b in {B/C, D/E, F/A}
\draw[->,line width=1pt,black] (\a) -- (\b);	
		
\filldraw [black] 
(A) circle (1.35pt)
(B) circle (1.35pt)
(C) circle (1.35pt)
(D) circle (1.35pt)
(E) circle (1.35pt)
(F) circle (1.35pt);
		
\end{tikzpicture}	
\begin{tikzpicture}[scale=0.73,>=Stealth]
		
\newcommand \x{1.5};
		
\draw[white] (-\x-0.5,-\x-0.5) rectangle (5*\x+0.5,\x+0.5);
		
\tikzset{->/.style={decoration={markings,mark=at position 1 with {\arrow[scale=0.83]{>}}},postaction={decorate}}}
		
\newcommand \q{1/1.17557050458};
		
\coordinate (A) at (\q*\x,0);
\coordinate (B) at ({cos(2*180 /5)*\q*\x},{-sin(2*180 /5)*\q*\x});
\coordinate (C) at ({-cos(180 /5)*\q*\x},{-sin(4*180 /5)*\q*\x});
\coordinate (D) at ({-cos(180 /5)*\q*\x},{sin(4*180 /5)*\q*\x});
\coordinate (E) at ({cos(2*180 /5)*\q*\x},{sin(2*180 /5)*\q*\x});
		
\coordinate (F) at (\q*\x+\x,0);
\coordinate (G) at (\q*\x+2*\x,0);
\coordinate (H) at (\q*\x+3*\x,0);
		
\coordinate (I) at ({\q*\x+((3+(sqrt(3) /2))*\x},0.5*\x);
\coordinate (J) at ({\q*\x+((3+(sqrt(3) /2))*\x},-0.5*\x);
		
\foreach \a/\b in {A/B, C/D, E/A, H/I, J/H}
\draw[->,line width=1pt,black!35!white] (\a) -- (\b);
		
\foreach \a/\b in {B/C, D/E, I/J}
\draw[->,line width=1pt,black] (\a) -- (\b);
		
\draw[line width=1pt,black] (A)--(F) (G)--(H);
\draw[{<[scale=0.83]}-{>[scale=0.83]},line width=1pt,black!35!white] (F)--(G);
		
\filldraw [black] 
(A) circle (1.35pt)
(B) circle (1.35pt)
(C) circle (1.35pt)
(D) circle (1.35pt)
(E) circle (1.35pt)
(F) circle (1.35pt)
(G) circle (1.35pt)
(H) circle (1.35pt)
(I) circle (1.35pt)
(J) circle (1.35pt);

\end{tikzpicture}
\caption{The two possible types of the (undirected) graph $G$,
  with $\cA$-edges shown in black, $\cB$-edges in gray.
  The arrows indicate the effect of the permutation $x_i \mapsto x_{\al(i)}$.}
\label{fig:g}
\end{figure}
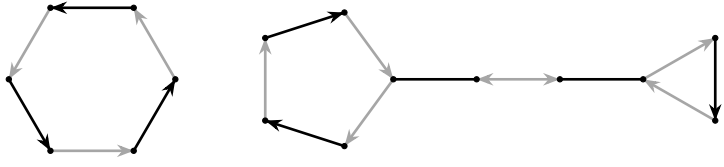

Now, in either case,
\begin{align*}
  \sum_{i \in I} d(x_i,x_{-i})
  &= \sum_{i \in I} \bigl( f(x_i) + f(x_{-i}) \bigr)\\
  &= \sum_{i \in I} \bigl( f(x_i) + f(x_{\al(i)}) \bigr)\\
  &= \sum_{i \in I} d(x_i,x_{\al(i)}) + k\bar s,  
\end{align*}
where $k$ is the number of indices $i \in I$ corresponding to
a $\cB$-edge $\{x_i,x_{\al(i)}\}$ of~$G$. We conclude that if $\bar s \le 0$,
then $X$ is $(n,0)$-hyperbolic. It remains to consider the case $\bar s > 0$.
Note that if $G$ is of the first type, then the cycle $G$ has at most
$|V| = 2(n+1)$ vertices, so $k \le n+1$. If $G$ is of the second type,
then $k$ equals the total number of $\cB$-edges in the two loops plus twice
the number of $\cB$-edges in the path from $v$ to $w$. This is equal to the
total number of $\cA$-edges in the two loops plus twice the number of
$\cA$-edges in the path from $v$ to $w$. Since each of the loops contains at
least one $\cA$-edge, which is counted only once, it follows that
$k \le |V| - 2 = 2n$. By Lemma~\ref{Lem:copy}, there is an isometric copy of
$\{\pm \bar s e_i\} \sub \li^{\,n+1}$ in $E(V)$, and hence also in $E(X)$.
Then $\bar s \le \del$ by assumption, so $k\bar s \le 2n\del$,
and therefore $X$ is $(n,n\del)$-hyperbolic. 
\end{proof}

Proposition~\ref{Pro:copy} may be viewed as a stable version of
Theorem~9 in~\cite{Dr} (p.~327), which follows as a corollary.

\begin{Cor} \label{Cor:dress}
A metric space $X$ is $(n,0)$-hyperbolic if and only if\/
$\dim_{\comb}(X) \le n$.
\end{Cor}

\begin{proof}
Let $n \ge 1$.
By Proposition~\ref{Pro:copy}, $X$ is $(n,0)$-hyperbolic if and only if
$E(X)$ contains no isometric copy of $\{\pm se_i\} \sub \li^{\,n+1}$ for
$s > 0$, and by Proposition~\ref{Pro:dim-comb} this holds if and only if
$\dim_\comb(X) \le n$. 
\end{proof}  

We now prove Theorem~\ref{Thm:nhyp-ih} stated in the introduction,
which subsumes some of the results obtained so far.

\begin{proof}[Proof of Theorem~\ref{Thm:nhyp-ih}]
If $X$ is $(n,\del)$-hyperbolic, then $E(X)$ is $(n,\del)$-hyperbolic
(Proposition~\ref{Pro:ih}) and has therefore asympotic rank $\le n$
(Proposition~\ref{Pro:asrk}). Evidently (3) implies~(4).
Lastly, suppose that~(4) holds. If $n = 0$, then $E(X)$ and $X$ are
bounded and so $X$ is $(0,\ast)$-hyperbolic. If $n \ge 1$, then $E(X)$
contains no isometric copy of $\{\pm s e_i\} \sub \li^{\,n+1}$ for
$s > 2r_0$ (Lemma~\ref{Lem:ball}), and hence
$X$ is $(n,2r_0n)$-hyperbolic (Proposition~\ref{Pro:copy}).
\end{proof}


\section{The slim simplex property and coarse injectivity}
\label{sect:ss-ci}

We turn to the slim simplex property~$(\SS_n)$ stated in the introduction.
By a {\em Euclidean $(n+1)$-simplex $\Del$} we mean the convex hull of $n+2$
points in $\R^{n+1}$ such that the interior of $\Del$ is non-empty, and a
{\em facet} of $\Del$ is the convex hull of $n+1$ of these vertices.
We restate Theorem~\ref{Thm:ss} in a slightly stronger form.
The proof uses Proposition~\ref{Pro:ih} and Proposition~\ref{Pro:asrk}
to apply a result from~\cite{GL}, which shows that~$(\SS_n)$ holds
uniformly for certain classes of proper metric spaces.

\begin{Thm} \label{Thm:sss}
Let $X$ be an $(n,\del)$-hyperbolic metric space.
Let $\Del$ be a Euclidean $(n+1)$-simplex, and let $\phi \colon \d\Del \to X$
be a map such that for some constants $\lam \ge 1$ and $c \ge 0$,
the restriction of $\phi$ to each facet of $\Del$ is a
$(\lam,c)$-quasi-isometric embedding. Then, for every facet $F$, 
the image $\phi(F)$ is contained in the closed $(1+c)D$-neighborhood of\/
$\phi(\overline{\d\Del \sm F})$ for some constant $D$ depending only on
$n,\del,\lam$ (and not on $X$).
\end{Thm}

\begin{proof}
If $n = 0$, then $\diam(X) \le \del$, and the result holds. Let now $n \ge 1$.
  
We consider $X$ as a subset of its injective hull $E(X)$ and write $d$
also for the metric of $E(X)$. First we approximate
$\phi \colon \d\Del \to X \sub E(X)$ by a piecewise Lipschitz map as follows.
Let $\beta$ denote the induced inner
metric on $\d\Del$. Since every shortest curve connecting two points in
$\d\Del$ meets each of the $n+2$ facets in at most one (possibly degenerate)
subsegment, $\phi$ is $(\lam,(n+2)c)$-Lipschitz with respect to~$\beta$. 
Let $Z \sub \d\Del$ be a maximal set subject to the condition that
$\beta(z,z') \ge (n+2)c\lam^{-1}$ whenever $z,z' \in Z$ are distinct.
For any such $z,z'$,
\[
  d(\phi(z),\phi(z')) \le \lam\,\beta(z,z') + (n+2)c \le 2\lam\,\beta(z,z').
\]
Since $E(X)$ is injective, $\phi|_Z$ extends to a $2\lam$-Lipschitz map
$\phi' \colon \d\Del \to E(X)$ with respect to $\beta$.
For every $x \in \d\Del$ there exists a $z \in Z$ with
$\beta(x,z) \le (n+2)c\lam^{-1}$, hence
\begin{align*}
  d(\phi'(x),\phi(x))
  &\le d(\phi'(x),\phi'(z)) + d(\phi(z),\phi(x)) \\
  &\le 2\lam\,\beta(x,z) + \lam\,\beta(x,z) + (n+2)c \\
  &\le 4(n+2)c.
\end{align*}
Furthermore, if $x,y$ are two points in the same facet, then
\[
  d(\phi'(x),\phi'(y)) \ge d(\phi(x),\phi(y)) - 8(n+2)c
  \ge \lam^{-1}\|x - y\| - c'
\]
for $c' := (8n + 17)c$, and $d(\phi'(x),\phi'(y)) \le 2\lam\,\beta(x,y)
= 2\lam\,\|x - y\|$.

By Proposition~\ref{Pro:ih}, $E(X)$ is $(n,\del)$-hyperbolic.
Since $\phi'(\d\Del) \sub E(X)$ is compact, so is its injective hull,
and hence there exists a compact injective subspace $Y \sub E(X)$ containing
$\phi'(\d\Del)$. We now apply Theorem~7.2 in~\cite{GL} for the class $\cX$
of all compact, injective, $(n,\del)$-hyperbolic spaces $Y$ (see also the
concluding remark in its proof for a simplification).
There are two assumptions on the class $\cX$ of metric spaces in
this theorem. The first is that all members of~$\cX$ satisfy
certain coning inequalities in dimensions $\le n$ with a uniform
constant. Since every injective metric space has a conical geodesic bicombing,
this holds with constant~$1$ (see Proposition~3.8 in~\cite{L},
Proposition~2.10 in~\cite{W05}, and Sect.~2.7 in~\cite{KL}).
The second assumption is that every sequence $(Y_k)_{k \in \N}$ in $\cX$
has asymptotic rank at most~$n$, which is satisfied by
Proposition~\ref{Pro:asrk}. The conclusion is that there is a constant
$D' = D'(\cX,n,\lam)$, hence depending only on $n,\del,\lam$, such that
for every facet $F$ of $\Del$, the image $\phi'(F)$ is contained in the
closed $(1+c')D'$-neighborhood of the union of the images of the remaining
facets. Since $\phi$ and $\phi'$ are uniformly close to each other,
this gives the result.
\end{proof}

Recall from the introduction that a metric space $X$ is 
{\em coarsely injective} if there exists a constant $c \ge 0$ such that
every $1$-Lipschitz map $\phi \colon A \to X$ defined on a subset of metric
space $B$ has a $(1,c)$-Lipschitz extension $\bar\phi \colon B \to X$.
To make the constant explicit, we say that $X$ is {\em $c$-coarsely injective}.
This property implies, more generally, that every
$(\lam,\eps)$-Lipschitz map $\phi \colon A \to X$ on $A \sub B$ has a
$(\lam,\eps + c)$-Lipschitz extension $\bar\phi \colon B \to X$,
because such a map $\phi$ is $1$-Lipschitz with respect
to the metric $d_{\lam,\eps}$ on $B$ satisfying $d_{\lam,\eps}(b,b') =
\lam\,d(b,b') + \eps$ for every pair of distinct points $b,b' \in B$.

The following result generalizes the well-known fact
that a metric space $X$ is injective if and only if $X$ is hyperconvex
(see~\cite{AP}). We call $X$ {\em coarsely hyperconvex} if, for some constant
$c \ge 0$, whenever $\{(x_s,r_s)\}_{s \in S}$ is a family in
$X \times \R$ satisfying $r_s + r_t \ge d(x_s,x_t)$ for all
pairs of indices $s,t \in S$, then $\bigcap_{s \in S}B(x_s,r_s + c) \ne \es$.
To make the constant explicit, we say that $X$ is
{\em $c$-coarsely hyperconvex}. For a geodesic metric space $X$, this
can be reformulated as the following {\em coarse Helly property} (compare
Sect.~3.3 in~\cite{CCGHO}):
any family $\{B(x_s,r_s)\}_{s \in S}$ of pairwise
intersecting closed balls in $X$ satisfies
$\bigcap_{s \in S}B(x_s,r_s+c) \ne \es$.

\begin{Pro} \label{Pro:c-inj}
For every metric space $X$, the following are equivalent:
\begin{enumerate}
\item
$X$ is coarsely injective;
\item
$X$ is coarsely hyperconvex;
\item
$E(X)$ is within finite distance from the image of the canonical embedding
$e \colon X \to E(X)$;
\item
$X$ is roughly isometric to an injective metric space $Y$.  
\end{enumerate}
\end{Pro}

In view of~(4) it is clear that all of these properties are preserved under
rough isometries. As the proof will show, all implications are quantitative.

For a geodesic Gromov hyperbolic space $X$, (1), (2), and~(3) were
established individually in \cite{L99}, \cite{CE}, and~\cite{L}, respectively. 
However, the equivalence of these properties was observed only recently;
see Proposition~3.12 in~\cite{CCGHO} for~(2) and~(3).

\begin{proof}
To show that (1) implies~(2), let $\{(x_s,r_s)\}_{s \in S} \sub X \times \R$
be a family such that $r_s + r_t \ge d(x_s,x_t)$ for all $s,t \in S$.
Consider the corresponding set $A := \{x_s: s \in S\}$ and put
$r(a) := \inf\{r_s + d(a,x_s): s \in S\}$ for all $a \in A$.
For $a,a' \in A$, the triangle inequalities
\[
  |r(a) - r(a')| \le d(a,a') \le r(a) + r(a')
\]
hold, thus there is a metric extension $B := A \cup \{b\}$ of $A$ with
$d(a,b) = r(a)$ for all $a  \in A$. Now if $X$ is $c$-coarsely injective,
then the inclusion map $A \to X$ extends to a $(1,c)$-Lipschitz map
on $B$, and the image point $y$ of $b$ satisfies 
$d(a,y) \le d(a,b) + c$ for all $a \in A$, hence $d(x_s,y) \le r_s + c$
for all $s \in S$.

We show that (2) implies~(3). Suppose that $X$ is $c$-coarsely hyperconvex,
and let $f \in E(X)$. Since $f(x) + f(x') \ge d(x,x')$ for all $x,x' \in X$,
it follows that there exists a point $y \in X$ with
$d(x,y) \le f(x) + c$ for all $x \in X$. Since $f$ is extremal,
$f(y) = \sup_{x \in X}(d(x,y) - f(x))$, thus (by~\eqref{eq:f-d})
$\|f - d_y\|_\infty = f(y) \le c$.

It is clear that if~(3) holds, then $e \colon X \to E(X)$
is a rough isometry between $X$ and~$E(X)$.

It remains to show that~(4) implies~(1).
Suppose that $i \colon X \to Y$ is a $(1,\eps)$-Lipschitz map into an
injective metric space $Y$, and $j \colon Y \to X$ is a $(1,\eps)$-Lipschitz
map such that $d(x,j \circ i(x)) \le \eps$ for all $x \in X$.
Let $\phi \colon A \to X$ be a $1$-Lipschitz map defined on $A \sub B$.
Then $i \circ \phi \colon A \to Y$ is $(1,\eps)$-Lipschitz and extends to a
$(1,\eps)$-Lipschitz map $\psi \colon B \to Y$, furthermore
$j \circ \psi \colon B \to X$ is $(1,2\eps)$-Lipschitz, and
\[
  d(\phi(a),j \circ \psi(a)) = d(\phi(a),j \circ i(\phi(a))) \le \eps
\]
for all $a \in A$. Hence, the map $\bar\phi \colon B \to X$ defined by
$\bar\phi(a) := \phi(a)$ for all $a \in A$ and
$\bar\phi(b) := j \circ \psi(b)$ for all $b \in B \sm A$ is a
$(1,3\eps)$-Lipschitz extension of $\phi$.
\end{proof}

We now prove our main result regarding coarsely injective spaces.

\begin{proof}[Proof of Theorem~\ref{Thm:cinj}]
Theorem~\ref{Thm:ss} shows that (1) implies~(2), and
Proposition~\ref{Pro:asrk} shows that (1) implies~(3) as well as~(4).

We show by contraposition that each of the conditions (2), (3), (4)
implies~(5). Suppose that there exist an $\eps > 0$ and
$(1,\eps)$-quasi-isometric embeddings of $B(0,k) \sub \li^{\,n+1}$ into $X$
for all integers $k \ge 1$.
Then one finds a Euclidean $(n+1)$-simplex $\Del$ and a sequence of 
maps $\phi_k \colon \d(k\Del) \to X$ violating~$(\SS_n)$. Furthermore,
the unit ball in $\li^{\,n+1}$ is an asymptotic subset of the constant
sequence $X_k = X$ and hence admits an isometric embedding into some
asymptotic cone $X_\om$ of $X$, thus $\asrk(X) \ge n+1$, and $X_\om$
fails to be $(n,0)$-hyperbolic.

For the proof of the implication (5) $\Rightarrow$ (1) and the last
assertion of the theorem, note that since $X$ is coarsely injective,
$E(X)$ is within finite distance from $e(X)$, so there exist a $c > 0$ and
a $(1,c)$-quasi-isometric embedding $E(X) \to X$.
Hence, if~(5) holds, then $E(X)$ cannot contain isometric copies of too
large balls in $\li^{\,n+1}$, and Theorem~\ref{Thm:nhyp-ih} then shows that
$X$ is $(n,\ast)$-hyperbolic. Similarly, if $X$ is proper and cocompact,
then there is an $n$ such that every set $V \sub X$ of distinct points at
mutual distance $\ge c$ and with diameter $\le 3c$ has less than
$2^{n+1}$ elements, thus $E(X)$ contains no isometric copy of
$\{-c,c\}^{\,n+1} \sub \li^{\,n+1}$. Then Theorem~\ref{Thm:nhyp-ih} shows that
$\asrk(E(X)) \le n$, and thus $\asrk(X) \le n$.
\end{proof}

We conclude with the proofs of the three corollaries stated in the
introduction.

\begin{proof}[Proof of Corollary~\ref{Cor:helly}]
Since $X$ is proper and cocompact, $n_4$ is finite and equal to the maximal
integer for which there exist isometric embeddings of
$B(0,r) \sub l_\infty^{\,n_4}$ into $X$ for all $r > 0$.
As an injective metric space, $X$ is isometric to $E(X)$, so
Theorem~\ref{Thm:nhyp-ih} shows that $n = n_1 = n_4$.
Since $G$ is quasi-isometric to $X$, Theorem~\ref{Thm:ss} shows further
that $n_2 \le n_1$, and evidently $n_4 \le n_3 \le n_2$.

For the last assertion, suppose to the contrary that $G$ has a free abelian
subgroup of rank $n+1$. By Proposition~3.8 in~\cite{L}, $X$ possesses an
equivariant conical geodesic bicombing, and it follows from Proposition~4.4
and Lemma~6.1 in~\cite{DL16} that there is an isometric embedding
of $\Z^{n+1}$ into $X$ with respect to the metric on $\Z^{n+1}$ induced by
some norm on $\R^{n+1}$. Thus there is a quasi-isometric embedding
of $\R^{n+1}$ into $X$. Alternatively, by the Algebraic Flat Torus Theorem for
semihyperbolic groups stated on p.~475 in~\cite{BH}, every monomorphism
of $\Z^{n+1}$ into $G$ is a quasi-isometric embedding.
\end{proof}

\begin{proof}[Proof of Corollary~\ref{Cor:hh}]
Since the rank $\nu$ of $X$ equals the `quasi-ball rank' $\bar\nu$,
for all $\lam,c$ there is a radius $r_0$ such that there is no
$(\lam,c)$-quasi-isometric embedding of $B(0,r) \sub \R^{\nu + 1}$ into $X$
for $r > r_0$. It follows that property~(5) of Theorem~\ref{Thm:cinj}
holds with $\nu$ in place of $n$ and with respect to the coarsely injective
metric $\rho$. Thus $(X,\rho)$ is $(\nu,\ast)$-hyperbolic, and $X$
satisfies~$(\SS_n)$ and has asymptotic rank at most $\nu$ with respect to
either $\rho$ or the original metric $d$. Furthermore, every metric space $X$
with quasi-ball rank $\bar\nu$ admits an asymptotic subset bi-Lipschitz
homeomorphic to the unit ball in $\R^{\bar\nu}$, so $\bar\nu \le \asrk(X)$
(in general the inequality may be strict; for example,
$X = \{k^2: k \in \N\} \sub \R$ satisfies $\bar\nu = 0$ and $\asrk(X) = 1$).
We conclude that $\nu = \asrk(X)$, and $\nu$ is the least integer
such that $(X,\rho)$ is $(\nu,\ast)$-hyperbolic.
\end{proof}

\begin{proof}[Proof of Corollary~\ref{Cor:symm}]
Every such Finsler metric $d$ is bi-Lipschitz equivalent to the Riemannian
metric, so $(X,d)$ has asymptotic rank $n$ and is therefore
$(n,\ast)$-hyperbolic by Theorem~\ref{Thm:cinj}.
\end{proof}

\subsection*{Acknowledgement}
We thank Marc Burger and Thomas Haettel for useful discussions and
for providing some references. 


\end{document}